\documentclass[12pt]{amsart}
\usepackage{amsmath,amsthm,amsfonts,amssymb,eucal}
\usepackage[notcite,notref]{showkeys}

\renewcommand {\a}{ \alpha }

\newcommand{\g}{\gamma}
\newcommand{\G}{\Gamma}

\newcommand{\varf}{\varphi}
\renewcommand{\d}{\delta}
\newcommand{\D}{\Delta}
\newcommand{\s}{\sigma}
\renewcommand{\l}{\lambda}
\renewcommand{\L}{\Lambda}
\newcommand{\z}{\zeta}

\newcommand{\p}{\partial}

\newcommand{\R}{ \mathbb R}
\newcommand{\C}{ \mathbb C}
\newcommand{\N}{ \mathbb N}

\newcommand{\CD}{\mathcal D}

\newcommand{\CF}{\mathcal F}

\newcommand{\CJ}{\mathcal J}
\newcommand{\CP}{\mathcal P}
\newcommand{\CR}{\mathcal R}
\newcommand{\CS}{\mathcal S}
\newcommand{\CT}{\mathcal T}

\newcommand{\CM}{\mathcal M}
\newcommand{\CN}{\mathcal N}

\newcommand {\GG}{\mathfrak G}
\newcommand {\GH}{\mathfrak H}

\newcommand {\GS}{\mathfrak S}

\newcommand {\gm}{\mathfrak m}

\newcommand {\bq}{\mathbf q}
\newcommand {\BA}{\mathbf A}

\newcommand {\BG}{\mathbf G}
\newcommand {\BH}{\mathbf H}

\newcommand {\BK}{\mathbf K}

\newcommand {\BQ}{\mathbf Q}

\newcommand{\CA}{\mathcal A}

\newcommand{\wt}{\widetilde}

 \DeclareMathOperator{\im}{Im}
\DeclareMathOperator{\re}{Re}

\newtheorem{thm}{Theorem}[section]

\newtheorem{lem}[thm]{Lemma}
\newtheorem{prop}[thm]{Proposition}

\theoremstyle{definition}
\newtheorem{defn}[thm]{Definition}

\theoremstyle{remark}

\numberwithin{equation}{section}

\newcommand{\thmref}[1]{Theorem~\ref{#1}}

\newcommand{\diag}{{\rm{diag}}}

\begin{document}

\title[spectrum of irreversible quantum graphs]
{Smilansky's model of irreversible quantum graphs, I: the absolutely
continuous spectrum}
\author[W.D. Evans]{W.D. Evans}
\address{School of Mathematics\\ Cardiff University\\
         23 Senghennydd Road\\ Cardiff CF24 4AG\\  UK}
\email{EvansWD@cardiff.ac.uk}
\author[M. Solomyak]{M. Solomyak}
\address{Department of Mathematics\\The Weizmann Institute of Science\\
Rehovot 76100\\Israel} \email{michail.solomyak@weizmann.ac.il}
\date{}

\begin{abstract}
In the model suggested by Smilansky \cite{SM} one studies an
operator describing the interaction between a quantum graph and a
system of $K$ one-dimensional oscillators attached at several
different points in the graph. The present paper is the first one
in which the case $K>1$ is investigated. For the sake of
simplicity we consider $K=2$, but our argument is of a general
character. In this first of two papers on the problem, we describe
the absolutely continuous spectrum. Our approach is based upon
scattering theory.
\end{abstract}

\maketitle
\section{Introduction}\label{int}

In the paper \cite{SM} U. Smilansky suggested a mathematical model
to which he gave the name ``Irreversible quantum graph''. In this
model one studies the interaction between a quantum graph and a
finite system of one-dimensional oscillators attached at several
different points in the graph. Recall that the term ``quantum
graph'' usually stands for a metric graph $\G$ equipped with a
self-adjoint differential operator acting on $L^2(\G)$; see the
survey paper \cite{Ku} and references therein. In our case this
operator will be the Laplacian $-\D$.

In Smilansky's model one initially deals with two independent
dynamical systems. One of the systems acts in $L^2(\G)$ and its
Hamiltonian is the Laplacian. Another system acts in the
space $L^2(\R^K),\ K\ge1$ and is generated by the  Hamiltonian
$H_{osc}={\sum}_{k=1}^K h_k$ where
\begin{equation*}\
h_k=\frac{\nu_k^2}{2}\left(-\frac{\p^2}{{\p q_k}^2}+
q_k^2\right),\qquad k=1,\ldots,K;
\end{equation*}
in \cite{SM} the oscillators are written in a slightly different
form; one form reduces to another by scaling. In what follows the
points in $\G$ are denoted by $x$ and the points in $\R^K$ by
$\bq=(q_1,\ldots,q_K)$.

Consider now the operator
\begin{equation}\label{1.1}
\BA_0=-\D\otimes I+I\otimes{H_{osc}}
\end{equation}
in the space $L^2(\G\times\R^K)$. It is defined by the differential expression
\begin{equation}\label{1.2}
\CA U=-\D_xU+\frac1{2}\sum_{k=1}^K\nu_k^2\bigl(-\frac{\p^2U}{{\p q_k}^2}
+q_k^2U\bigl)
\end{equation}
and is self-adjoint on the natural domain. The terms in \eqref{1.1} do not
interact with each other.

Interaction is introduced with the help of a system of
``matching conditions'' on the derivative $U'_x$ at some points
$o_1,\ldots,o_K\in\G$. One says that the $k$-th oscillator
is attached to the graph at the point $o_k$.
The condition at the point $o_k$ is
\begin{equation}\label{1.4}
[U'_x](o_k,\bq)=\a_kq_kU(o_k,\bq),\qquad k=1,\ldots,K,
\end{equation}
where $[f'_x](.)$ stands for the expression appearing in the
Kirchhoff condition, well known in the theory of electric networks.
When $\G=\R$ (which is the only case we deal with in the main
body of the paper), $[f'_x](.)$ is the jump of the derivative,
\begin{equation}\label{1.5}
[f'_x](o)=f'_x(o+)-f'_x(o-).
\end{equation}
The real parameter $\a_k$ in \eqref{1.4} expresses the strength of
interaction between the quantum graph and the oscillator $h_k$.
The case
$\a_1=\ldots=\a_K=0$ corresponds to the operator $\BA_0$ as in \eqref{1.1}.

Sometimes we shall denote by $\a,\nu$ the multi-dimensional
parameters $\a=\{\a_1,\ldots,\a_K\}$,
$\nu=\{\nu_1,\ldots,\nu_K\}$. Let
$\BA_{\a;\nu}=\BA_{\a_1,\ldots,\a_K;\nu_1,\ldots,\nu_K}$ stand for
the operator defined by the differential expression \eqref{1.2}
and the conditions \eqref{1.4}. Usually, the values of $\nu_k$ are
fixed and we exclude them from the notation. On the other hand, we
use the notation $\BA_{\G;\a;\nu}$ for this operator when it is
necessary to reflect its dependence on the graph.

The problem to be considered is the description of the spectrum of
the dynamical system generated by the Hamiltonian $\BA_{\a;\nu}$.
More specifically, it is to construct the self-adjoint realization
of $\BA_{\a;\nu}$ as an operator in the Hilbert space
$L^2(\G\times\R^K)$ and to describe its spectrum. \vskip0.2cm Up
until now, the problem has only been investigated for the simplest
case $K=1$. The first results were obtained in the paper \cite{SM}
by Smilansky. Then a detailed study of the problem was carried out
in the papers \cite{S1}, \cite{S2} and \cite{NS}. In \cite{S3},
along with some new results, a detailed survey of the current
state of the problem is given.

\vskip0.2cm On first sight, the
problem might seem amenable to the perturbation theory of
quadratic forms. Indeed, the spectrum $\s(\BA_0)$ can be easily
described by separation of variables and the perturbation in the
quadratic form, which appears when passing from $\BA_0$ to
$\BA_\a$ with $\a\neq0$, seems not to be too strong. However, this
is not so: this perturbation turns out to be only form-bounded but
not form-compact, which makes it impossible to apply the standard
techniques. So, the problem requires certain specific tools which
were developed in \cite{S1} -- \cite{S3} and \cite{NS}. The most
important of these tools is the systematic use of Jacobi matrices.
\vskip0.2cm It was found in the above mentioned papers on the
one-oscillator problem that the character of the spectrum strongly
depends on the size of $\a$: there exists some $\a^*>0$ such that
the absolutely continuous spectrum $\s_{a.c.}(\BA_\a)$ coincides
with $\s_{a.c.}(\BA_0)$ if $|\a|<\a^*$ (in particular, it is
absent if the graph is compact) and fills the whole of $\R$ if
$|\a|>\a^*$. The dependence of the structure of the point spectrum
$\s_{p}(\BA_\a)$ on $\a$ is also well understood. \vskip0.2cm

This is the first of two papers on the problem for $K>1$ and in it
we study the absolutely continuous spectrum; in our other paper
\cite{ES} the point spectrum is investigated. This division is
natural, since the technical tools used in each part are
different. We address the simplest situation, when $\G=\R$ and
$K=2$, but our argument is of a rather general character and we
firmly believe that it applies to a wide class of graphs and to
any $K$. However, in the general case, the calculations become
more complicated and this obscures the main features of the
argument.

\vskip0.2cm We first describe informally the main ideas
lying behind our approach.\vskip0.2cm

The effect of adding one more oscillator to a system with $K$
oscillators is twofold. Firstly, the total dimension of the set
$\G\times\R^{K}$ increases by one which certainly affects the
spectrum. Secondly, there is some effect coming from the
additional matching condition \eqref{1.4} at the point $o_{K+1}$.
This second effect disappears if we take $\a_{K+1}=0$. Indeed,
then the variable $q_{K+1}$ can be separated and the operator
decomposes into the orthogonal sum of simpler operators. More
exactly, denote by $\wt{\BA}$ the operator which corresponds to
the configuration with the $(K+1)$-th oscillator removed,
\begin{equation*}
\wt{\BA}=\BA_{\a_1,\ldots,\a_K;\nu_1,\ldots,\nu_K}.
\end{equation*}
Then it is easy to see that
\begin{equation}\label{1.6}
\BA_{\a_1,\ldots,\a_K,0;\nu_1,\ldots,\nu_K,\nu_{K+1}}
={\sum_{n\in\N_0}}^\oplus(\wt{\BA}+\nu^2_{K+1}(n+1/2)).
\end{equation}
This orthogonal decomposition yields the complete description of
the spectrum of the operator on the left-hand side, provided that
the spectrum of $\wt{\BA}$ is known.

The key observation which allows one to solve the general problem
is that the interaction between the oscillators attached at
different points is weak. For $K=2$ this observation leads to the
conclusion that the study of $\s(\BA_{\a_1,\a_2;\nu_1,\nu_2})$ can
be reduced to the same problem for the operators
$\BA_{\a_1,0;\nu_1,\nu_2}$ and $\BA_{0,\a_2;\nu_1,\nu_2}$. Due to
the equality \eqref{1.6} this reduces the problem to the study of
the spectra of two operators, $\BA_{\a_1;\nu_1}$ and
$\BA_{\a_2;\nu_2}$, each corresponding to the case of only one
oscillator. Since the latter case is already well understood, we
obtain the desired results for our more complicated case.
\vskip0.2cm An accurate realization of this idea is different for
the point spectrum and for the absolutely continuous spectrum. In
the present paper we concentrate on the absolutely continuous
spectrum. Here an important correction to the above
scheme is necessary: the
study of $\s_{a.c.}(\BA_{\a_1,\a_2;\nu_1,\nu_2})$ does not reduce
to the study of $\s_{a.c.}(\BA_{\a_1;\nu_1})$ and
$\s_{a.c.}(\BA_{\a_1;\nu_1})$ for the same graph $\G$. Rather, we
have to divide $\G$ into two parts, $\G=\G_1\cup\G_2$ in such a
way that $o_j\in\G_j$ and $o_j\notin\G_{3-j}$. Then
$\s_{a.c.}(\BA_{\G;\a_1,\a_2;\nu_1,\nu_2})$ can be expressed in
terms of $\s_{a.c.}(\BA_{\G_j;\a_j;\nu_j})$, $j=1,2$.
The paper \cite{ES} is
devoted to the study of the point spectrum.
There such a partition of $\G$ is unnecessary.

\bigskip

We use the following notation. We write $\N_0$ for the set
$\{0,1,\ldots\}$.
The diagonal operator in an appropriate
$\ell^2$-space, with the diagonal elements $a_0,a_1,\ldots$, is
denoted by $\diag\{a_n\}$. We apply similar notation for the
block-diagonal operators. The notation $\CJ(\{a_n\},\{b_n\})$
stands for the Jacobi matrix whose non-zero entries are
$j_{n,n}=a_n$ and $j_{n,n+1}=j_{n+1,n}=b_n$. If $\BA$ is a
self-adjoint operator in a Hilbert space, then $\s(\BA),
\s_{a.c}(\BA), \s_p(\BA)$ stand for its spectrum, absolutely
continuous (a.c.) spectrum and point spectrum respectively. We use
the symbol $\gm_{a.c}(\l;\BA)$ for the multiplicity function of
the a.c. spectrum. The symbol $\GS_1$ stands for the trace class
of compact operators.

Other necessary notations are introduced in the course of the
presentation.

\section{Statement of the problem. Results}\label{st}
\subsection{The operator $\BA_\a$.}\label{st1}
As was mentioned in the introduction, we present our argument for
the graph $\G=\R$ and $K=2$.
We choose the points $o_1=1,o_2=-1$ and denote the coordinates in $\R^2$
by $q_+,q_-$ and
the parameters by $\a=\{\a_+,\a_-\}$, $\nu=\{\nu_+,\nu_-\}$. The
Laplacian on $\G$ is just the operator $-d^2/dx^2$ with the
Sobolev space $H^2(\R)$ as the operator domain. The operator $\BA_{\a,\nu}$
acts in the Hilbert space $\GH=L^2(\R^3)$ and
is defined by the differential expression
\begin{equation}\label{st.1}
\CA U=\CA_\nu U=-U''_{x^2}+\frac{\nu_+^2}{2}(-U''_{q_+^2}+q_+^2U)+
\frac{\nu_-^2}{2}(-U''_{q_-^2}+q_-^2 U)
\end{equation}
and the matching conditions (cf. \eqref{1.5})
\begin{equation}\label{st2}
[U'_x](\pm1,q_+,q_-)=\a_\pm q_{\pm}U(\pm1,q_+,q_-).
\end{equation}
So, in the notation of the introduction, we are dealing with the
operator
\[ \BA_{\R;\a_+,\a_-;\nu_+,\nu_-}.\]
However, as a rule we use the shortened notation $\BA_\a$. Note
that the replacement $\a_\pm\mapsto -\a_\pm$ corresponds to the change
of variables $q_\pm\mapsto -q_\pm$ which does not affect the spectrum.
For this reason, we discuss only
$\a_\pm\ge 0$.

The structure of the differential expression $\CA$ makes it
natural to decompose the function $U$ in a double series in
terms of the normalized Hermite functions $\chi_n$, namely
\begin{equation}\label{sa.0}
U(x,q_+,q_-) = \sum_{m,n \in \N_0}u_{m,n}(x)\chi_m (q_+) \chi_n(q_-),
\end{equation}
which is hereafter represented by $U \sim \{u_{m,n}\}$. The
mapping $U\mapsto\{u_{m,n}\}$ is an isometry of the space
$L^2(\R^3)$ onto the Hilbert space $\GH=\ell^2(\N_0^2;L^2(\R))$. We
evidently have $\CA U\sim\{L_{m,n}u_{m,n}\}$ where
\begin{equation}\label{sa.1}
(L_{m,n}u)(x)=-u''(x)+r_{m,n}u(x),\qquad x\neq\pm1;
\end{equation}
\begin{equation}\label{sa.2}
r_{m,n}=\nu_+^2(m+1/2)+\nu_-^2(n+1/2),\qquad m,n\in\N_0.
\end{equation}
The conditions at $x=\pm1$ reduce to
\begin{equation}\label{sa.3}\begin{split}
[u'_{m,n}](1)&=\frac{\a_+}{\sqrt 2}\left(\sqrt{m+1}u_{m+1,n}(1)+
\sqrt{m}u_{m-1,n}(1)
\right);\\
[u'_{m,n}](-1)&=\frac{\a_-}{\sqrt 2}\left(\sqrt{n+1}u_{m,n+1}(-1)+
\sqrt{n}u_{m,n-1}(-1)\right).\end{split}
\end{equation}
To derive the conditions \eqref{sa.3} from \eqref{st2}, one uses
the recurrency equation for the functions $\chi_n$,
\begin{equation*}
\sqrt{n+1}\chi_{n+1}(q)-\sqrt{2}q\chi_{n}(q)+\sqrt{n}\chi_{n-1}(q)=0.
\end{equation*}

\vskip0.2cm
\subsection{Operator $\BA_0$.}\label{op0}
The operator $\BA_0:=\BA_{0,0;\nu_+,\nu_-}$
admits separation of variables and we get
\begin{equation}\label{add1}
\BA_0=\sum_{m,n}^\oplus(\BH_0+r_{m,n}),
\end{equation}
where $\BH_0$ is the self-adjoint operator $-d^2/dx^2$ in
$L^2(\R)$. This leads to the complete description of the spectrum
$\s(\BA_0)$, namely, that it is purely a.c. and fills the
half-line $[r_{0,0},\infty)=[(\nu_+^2+\nu_-^2)/2,\infty)$. The
expression for the multiplicity function $\gm_{a.c.}(\l;\BA_0)$
immediately follows from \eqref{add1}, but is omitted.

\subsection{Domain of  $\BA_\a$.}\label{opa}
It is convenient to describe the domain of the self-adjoint
realization of the operator $\BA_\a$ in terms of the decomposition
\eqref{sa.0}. Define the set $\CD_\a$ as follows.
\begin{defn}\label{d}
An element $U\sim\{u_{m,n}\}$ lies in $\CD_\a$ if and only if

\noindent 1. $u_{m,n}\in H^1(\R)$ for all $m,n$.

\noindent 2. For all $m,n$ the restriction of $u_{m,n}$ to each
interval $(-\infty,-1),\ (-1,1)$, $(1,\infty)$ lies in $H^2$ and
moreover,
\begin{equation*}
\sum_{m,n}\int_\R |L_{m,n}u_{m,n}|^2dx<\infty.
\end{equation*}
3. The conditions \eqref{sa.3} are satisfied.
\end{defn}
Along with the set $\CD_\a$, define its subset
\begin{equation*}
\CD_\a^\bullet=\left\{ U\in\CD_\a: U\sim\{u_{m,n}\} \ \rm{finite}\right\}
\end{equation*}
where by {\it finite} we mean that the sequence has only a finite number of
non-zero components.
Denote by $\BA^\bullet_\a$ the operator in $\GH=L^2(\R^3)$, defined by
the system \eqref{sa.1} on the domain $\CD^\bullet_\a$,

\begin{lem}\label{l1}
The operator $\BA_\a^\bullet$ is symmetric in $\GH$. Its adjoint coincides
with the operator $\BA_\a$ considered on the domain $\CD_\a${\rm :}
\[(\BA_\a^\bullet)^*=\BA_\a.\]
\end{lem}
The proof is a straightforward modification of that for (5.2) in \cite{NS}.
\vskip0.2cm

\begin{thm}\label{sa.thm} For any $\a_+,\a_-\ge 0$ the operator
$\BA_\a$ is self-adjoint.
\end{thm}
The proof is given in section \ref{pr}. \thmref{sa.thm} and Lemma \ref{l1}
show that $\BA_\a$ is the unique natural self-adjoint realization of the
operator, defined by the differential expression \eqref{st.1} and the
matching conditions \eqref{st2}.

\subsection{Absolutely continuous spectrum of
the operator $\BA_\a$.}\label{circ}
Below we construct an
operator $\BA^\circ_\a$ whose a.c. spectrum admits a complete description.
Then we show
that the a.c. spectra of both operators $\BA_\a$ and $\BA^\circ_\a$
coincide, including the
multiplicities.

As a first step, let us consider two operators, $\BA^+_{\a_+}$ and
$\BA^-_{\a_-}$. The operator $\BA^+_{\a_+}$, say, acts in
the space $L^2(\R_+\times\R^2)$
and is defined by
the differential expression \eqref{st.1}, the matching condition \eqref{st2}
at the point $o_1=1$ and the Dirichlet condition $U(0,q_+,q_-)=0$.
The definition of $\BA^-_{\a_-}$ is similar, with $\R_+$ replaced by $\R_-$
and the point $o_1=1$ by $o_2=-1$.
By separation of variables, the operators $\BA^\pm_{\a_\pm}$ can be identified
with the orthogonal sum of simpler operators:

\begin{equation}\label{one.2}\begin{split}
\BA^+_{\a_+}={\sum_{n\in\N_0}}^\oplus
\left(\BA_{\R_+;\a_+;\nu_+}+\nu_-^2(n+1/2)\right)&,\\
\BA^-_{\a_-}={\sum_{m\in\N_0}}^\oplus
\left(\BA_{\R_-;\a_-;\nu_-}+\nu_+^2(m+1/2)\right)&.\end{split}
\end{equation}
Hence, both operators are self-adjoint. The direct sum
\begin{equation}\label{one.3}
\BA^\circ_{\a}=\BA^\circ_{\a_+,\a_-;\nu_+,\nu_-}:=\BA^+_{\a_+}\oplus\BA^-_{\a_-}
\end{equation}
is a self-adjoint operator in the original Hilbert space $\GH$.

\vskip0.2cm The following theorem is the main result of the paper.
Its formulation involves the notion of wave operator, which is one
of the basic notions in mathematical scattering theory; see e.g. \cite{K},
\cite{RS} or \cite{Y}.

\begin{thm}\label{one.t1} For each of the pairs $(\BA_\a,\ \BA^\circ_\a),
(\BA^\circ_\a,\ \BA_\a)$, there exist complete isometric wave
operators. In particular, the absolutely continuous parts of
$\BA_{\a}$ and $\BA^\circ_\a$ are unitarily equivalent.
\end{thm}

\thmref{one.t1} and the formulae \eqref{one.2}, \eqref{one.3}
reduce the study of
$\s_{a.c.}(\BA_\a)$ to the similar problem for the case of only
one oscillator. The latter problem was solved in \cite{NS} and
\cite{S3}. The next statement collects, for the particular case we
need, the results of section 3 in \cite{S3}; see also Theorem 5.1
and remarks in section 9 of \cite{NS}. In both papers it was
assumed that $\nu=1$, and we arrive at the formulation below via
scaling. By default, we take $\gm_{a.c.}(\l;\BA)=0$ if
$\l\notin\s_{a.c.}(\BA)$.

\begin{prop}\label{prop}{\rm{(The case of one oscillator.)}}
Let $\G=\R_+$ and $o=1$, or $\G=\R_-$ and $o=-1$. Then

\noindent {\rm 1)} \hskip 3.5cm
$\s_{a.c.}(\BA_{0;\nu})=[\nu^2/2,\infty)$;
\[\gm_{a.c.}(\l;\BA_{0;\nu})=n\ {\rm{for}}\ -\nu^2/2\le\l-\nu^2 n<\nu^2/2,\ n\in\N
\rm;\]

\noindent {\rm2)} if $0<\a < \nu\sqrt2 $, then \hskip 3.5cm
\[\s_{a.c.}(\BA_{\a;\nu})=\s_{a.c.}(\BA_{0;\nu}) =
[\nu^2/2,\infty);\]
 \[\gm_{a.c.}(\l;\BA_{\a;\nu})=\gm_{a.c.}(\l;\BA_{0;\nu});\]

\noindent {\rm3)} if $\a=\nu\sqrt2 $, then
\[\s_{a.c.}(\BA_{\a;\nu})=[0,\infty)\rm;\qquad
\gm_{a.c.}(\l;\BA_{\a;\nu})=\gm_{a.c.}(\l;\BA_{0;\nu})+1,\ \forall\l\ge 0;\]

\noindent {\rm4)} if $\nu\sqrt2 <\a<\infty$, then
\[\s_{a.c.}(\BA_{\a;\nu})=\R;\qquad
\gm_{a.c.}(\l;\BA_{\a;\nu})=\gm_{a.c.}(\l;\BA_{0;\nu})+1,\ \forall\l\in\R.\]
\end{prop}

Now we are in a position to present the final formula for the
function $\gm_{a.c.}(\l;\BA^\circ_{\a})$, and thus for our
original operator $\BA_\a$.
\begin{equation}\label{add3}\begin{split}
\gm_{a.c.}(\l;\BA_{\a})=&
\sum_{n\in\N_0}\gm_{a.c.}(\l-\nu_-^2(n+1/2);\BA_{\R_+;\a_+;\nu_+})\\
+&\sum_{m\in\N_0}\gm_{a.c.}(\l-\nu_+^2(m+1/2);\BA_{\R_-;\a_-;\nu_-}).\end{split}
\end{equation}

Combining the equality \eqref{add3} with Proposition \ref{prop},
we obtain the following description of the a.c. spectrum of the
operator $\BA_{\a}$ for any $\a_+,\a_-\ge 0$.

\begin{thm}\label{one.t2a} Let $\BA_{\a}= \BA_{\a;\nu}$ be
the self-adjoint operator defined by the
differential expression \eqref{st.1} on the operator domain $\CD_\a$.

\noindent {\rm 1)} If $\a_{\pm}/\nu{\pm}<\sqrt 2$, then
\[ \s_{a.c.}(\BA_{\a})=[r_{0,0},\infty)=[(\nu_+^2+\nu_-^2)/2,\infty).\]

\noindent {\rm2)} Let $\a_+/\nu_+ =\sqrt 2 $ and $\a_-/\nu_-
<\sqrt 2$, or $\a_-/\nu_- =\sqrt 2$ and $\a_+/\nu_+ <\sqrt 2$.
Then
\[ \s_{a.c.}(\BA_{\a})=[\nu_-^2/2,\infty)\qquad {\rm{or}}\qquad
 \s_{a.c.}(\BA_{\a})=[\nu_+^2/2,\infty)\]
respectively.

\noindent {\rm3)} Let $\a_+/\nu_+ =\a_-/\nu_- =\sqrt 2$, then
\[ \s_{a.c.}(\BA_{\a})=[0,\infty).\]

In all the cases {\rm 1 -- 3} the multiplicity function $\gm_{a.c.}(\l;\BA_{\a})$,
given by the equality \eqref{add3}, is finite for
all $\l\in\s_{a.c.}(\BA_{\a})$.
\vskip0.2cm
\noindent {\rm4)} Let $max(\a_+/\nu_+,\a_-/\nu_-)>\sqrt2$. Then
\[\s_{a.c.}(\BA_{\a})=\R,\qquad \gm_{a.c.}(\l;\BA_{\a})\equiv\infty.\]
\end{thm}
In connection with this theorem, we would like to emphasize that
the existence of the wave operators established in \thmref{one.t1}
gives much more information about the operator $\BA_\a$ than just
the description of its a.c. spectrum.

For the proof of  \thmref{one.t1} we use the following classical result
due to Kato, see Theorem 6.5.1 and Remark 6.5.2 in \cite{Y}.

\begin{prop}\label{kato} Let $\BA, \BA^\circ$ be self-adjoint operators in a
Hilbert space. Suppose that for some natural number $p$ the inclusion
\[ (\BA^\circ-\L)^{-p}-  (\BA^\circ-\L)^{-p}\in\GS_1\]
is satisfied for all non-real $\L\in\C$. Then the complete
isometric wave operators exist for both pairs $\BA, \BA^\circ$ and $\BA^\circ,\BA$.
\end{prop}

In our case the conditions of Proposition \ref{kato} turn out to be fulfilled with $p=3$.
This is the result of the following statement whose proof is our main technical goal
in this paper.
\begin{thm}\label{one.t2} For any non-real $\L\in\C$ one has
\begin{equation}\label{goal}
(\BA^\circ_\a-\L)^{-3} - (\BA_\a-\L)^{-3}\in \GS_1.
\end{equation}
\end{thm}
\thmref{one.t1} is a direct consequence of \thmref{one.t2}.

The proof of \thmref{one.t2} is rather long and requires some
preparatory work.

\section{Auxiliary material}\label{au}
In this section we present some elementary technical material
concerning the equations
\begin{gather}
-u''+\z^2u=f,\label{res.1g}\\
-v''+\z^2 v=0.\label{res.1go}
\end{gather}
where $\z=\g+i\d$ is a complex parameter. We need this material for the proofs
of both our main technical results, \thmref{sa.thm} and \thmref{one.t2}.

We assume that $\g>0$, and are mainly interested in estimates
which are uniform with respect to $\z$.

\subsection{Homogeneous equation.}\label{he}
Let $\CF_\z$ be the two-dimensional
space of functions on $\R$ which are continuous,
vanish as $|x|\to\infty$, and for $x\neq\pm1$ satisfy the equation \eqref{res.1go}.
We choose the following basis $\varf_\z^+,\varf_\z^-$ in $\CF_\z$:
\begin{equation*}
\varf_\z^+(x)=\begin{cases} 0,& x<-1,\\ \frac{\sinh\,\z(x+1)}{\sinh\,2\z},&
-1\le x\le 1;\\ e^{-\z(x-1)},& x>1;\end{cases}\qquad
\varf_\z^-(x)=\varf_\z^+(-x).
\end{equation*}
Then
\begin{equation}\label{1}
\varf_\z^+(1)=\varf_\z^-(-1)=1,\qquad \varf_\z^+(-1)=\varf_\z^-(1)=0.
\end{equation}
Just for this reason this basis is more convenient than the ``natural'' basis
consisting of the functions $e^{-\z|x\pm1|}$. Note also that
\begin{equation*}
\left[(\varf_\z^\pm)'\right](\pm1)=-\frac{2\z}{1-e^{-4\z}},\qquad
\left[(\varf_\z^\pm)'\right](\mp1)= \frac{2\z e^{-2\z}}
{1-e^{-4\z}}
\end{equation*}
and hence, for all $v\in\CF_\z$,
\begin{equation}\label{1c}\begin{split}
[v'](1)=&-\frac{2\z}{1-e^{-4\z}}\left(v(1)-e^{-2\z}v(-1)\right),\\
[v'](-1)=&-\frac{2\z}{1-e^{-4\z}}\left(v(-1)-e^{-2\z}v(1)\right).\end{split}
\end{equation}

A standard calculation shows that for the norm and scalar product
in $L^2(\R)$,
\begin{gather}
\|\varf_\z^+\|^2=\|\varf_\z^-\|^2=\frac1{2\g}+\frac{\g^{-1}\sinh 4\g-
\d^{-1}\sin 4\d}{4(\sinh^2 2\g+\sin^2 2\d)}=\g^{-1}+o(e^{-4\g}),\label{1f}\\
(\varf_\z^+,\varf_\z^-)=O(e^{-2\g}),\qquad \g\to\infty.\notag
\end{gather}
This shows that for $\g$ large the chosen basis is ``almost
orthogonal''. It follows that the two-sided estimate
\begin{equation}\label{1x}
c_0^{-1}\g\|v\|^2\le|C_+|^2+|C_-|^2\le c_0\g\|v\|^2,
\qquad v=C_+\varf_\z^+ +C_-\varf_\z^-\in\CF_\z
\end{equation}
with some $c_0>1$
is satisfied uniformly in any half-plane $\g=\re\z\ge\g_0>0$.

\vskip0.2cm
Now we turn to the subspace $\CF^\circ_\z$ formed by the functions
$v\in\CF_\z$, satisfying an additional condition $v(0)=0$. The functions
\begin{equation*}
\varf_\z^{\circ,+}(x)=\begin{cases}0,&x<0,\\ \frac{\sinh\z x}{\sinh\z},
&0\le x\le 1,\\
e^{-\z(x-1)},& x>1;\end{cases}\qquad
\varf_\z^{\circ,-}(x)=\varf_\z^{\circ,+}(-x)
\end{equation*}
form a natural basis in $\CF^\circ_\z$. For the functions $\varf_\z^{\circ,\pm}$
the equalities \eqref{1} are satisfied, and instead of \eqref{1c} we have
\begin{equation}\label{1ac}
[v'](\pm1)=-\frac{2\z}{1-e^{-2\z}}v(\pm1),\qquad\forall
v\in\CF^\circ_\z.
\end{equation}
Similarly to \eqref{1f}, we find that
\begin{equation*}
\|\varf_\z^{\circ,\pm}\|^2=\g^{-1}+O(e^{-2\g}),\qquad
\left(\varf_\z^{\circ,+},\varf_\z^{\circ,-}\right)=O(e^{-\g}),
\qquad \g\to\infty. \end{equation*} As a consequence, we conclude
that an analogue of \eqref{1x}, with the functions
$\varf_\z^{\pm}$ replaced by $\varf_\z^{\circ,\pm}$, is valid for
$v\in\CF^\circ_\z$. \vskip0.2cm

A straightforward calculation shows also that
\begin{equation}\label{2c}
\|\varf_\z^{\circ,\pm}-\varf_\z^{\pm}\|=O(e^{-\g}),\qquad \g\to\infty.
\end{equation}

\subsection{Non-homogeneous equation.}\label{nhom}
Here we discuss the equation \eqref{res.1g} without the matching
conditions at $x=\pm1$ or, equivalently, under the conditions of
the type \eqref{sa.3} with $\a=0$. Then the solution is given by
\begin{equation}\label{res.1}
u_\z(x)=(2\z)^{-1}\int_\R e^{-\z|x-t|}f(t)dt.
\end{equation}

The solution of the same equation \eqref{res.1g} subject to the condition $u(0)=0$
is
\begin{equation}\label{res.1b}
u^\circ_{\z}(x)=\begin{cases}
(2\z)^{-1}\int_{\R_+}\left(e^{-\z|x-t|}-e^{-\z(x+t)}\right)
f(t)dt,\qquad x>0;\\
(2\z)^{-1}\int_{\R_-}\left(e^{-\z|x-t|}-e^{\z(x+t)}\right)
f(t)dt,\qquad x<0.\end{cases}
\end{equation}

\vskip0.2cm

The difference $u^\circ_{\z}-u_{\z}$ is given by a rank one operator,
\begin{equation}\label{res.1f}
u^\circ_{\z}(x)-u_{\z}(x)=-(2\z)^{-1}g_\z(x)\int_{\R}f(t)g_\z(t)dt,\qquad
g_\z(x)=e^{-\z|x|}.
\end{equation}
Note that $\|g_\z\|^2=\g^{-1}$.

\subsection{Dependence on the additional parameters.}\label{add}
We are particularly interested in the case when $\z$ depends on
two parameters $r\in\R$  and $\L\in\C$, where $r\ge r_0>0,\
\L\notin\R_+$:
\begin{equation}\label{2}
\z=\z_r(\L):=\g_r(\L)+i\d_r(\L)=\sqrt{r-\L}.
\end{equation}
We select the branch of the square root in \eqref{2} to have
\begin{equation*}
\re\z_r(\L)>0,\qquad \im\L\cdot\im\z_r(\L)\le0.
\end{equation*}

For $\L$ fixed all the points $\z_r(\L)$ lie in some half-plane
$\re \z_r(\L)\ge\g_0(\L)>0$, hence \eqref{1x} is satisfied. It is
clear that $\g_r(\L)\sim r^{1/2}$ as $r\to\infty$. Therefore,
for any $\L\notin\R$
there exists a constant $c_1=c_1(\L)>1$ such that
\begin{equation}\label{2r}\begin{split}
c_1^{-1}r^{1/2} \|v\|^2\le|C_+|^2+|C_-|^2&\le c_1 r^{1/2}\|v\|^2, \\
\forall v=C_+\varf_\z^+ +C_-\varf_\z^-\in\CF_\z,&\qquad \z=\z_r(\L).\end{split}
\end{equation}

\section{Self-adjointness: proof of \thmref{sa.thm}}\label{pr}
According to the general theory of self-adjoint operators,
we must show that the equation $\BA_\a V=\L V$ has only the trivial solution
for some (and then all) $\L\in\C_\pm$.
To simplify our notation, we shall denote
\begin{equation}\label{sa.4}
\z_{m,n}(\L)=\z_{r_{m,n}}(\L)=(r_{m,n}-\L)^{1/2},\qquad
\varf_{m,n}^\pm(x;\L)=\varf_{\z_{m,n}(\L)}^\pm(x).
\end{equation}
If $V\sim\{v_{m,n}\}$,
then each function $v_{m,n}$ can be written as
\begin{equation}\label{sa.4y}
v_{m,n}(x)=r_{m,n}^{1/4}\left(C_{m,n}^+\varf_{m,n}^+(x;\L)+
C_{m,n}^-\varf_{m,n}^-(x;\L)\right),
\end{equation}
with some coefficients $C_{m,n}^\pm$. We have inserted the factor
$r_{m,n}^{1/4}$ in order that (cf. \eqref{2r})
\begin{equation}\label{sa.4a}
\{C_{m,n}^+,C_{m,n}^-\}\in\ell^2\,\Longleftrightarrow \{V\in\GH\}.
\end{equation}
The matching conditions \eqref{sa.3} at $x=\pm1$ yield an infinite
system of homogeneous linear equations for the unknown
coefficients $C_{m,n}^\pm$. Below we set $\mu_\pm=\sqrt2 /\a_\pm$.
Taking \eqref{1c} into account, we get from the condition at
$x=1$:
\begin{gather}
r_{m+1,n}^{1/4}(m+1)^{1/2}C_{m+1,n}^+ +
\frac{2\mu_+\z_{m,n}(\L)r_{m,n}^{1/4}}{1-e^{-4\z_{m,n}(\L)}}
\left(C_{m,n}^+ -C_{m,n}^-e^{-2\z_{m,n}(\L)}\right)\notag\\
+r_{m-1,n}^{1/4}m^{1/2}C_{m-1,n}^+=0.
\label{sa.4b}\end{gather}
It is convenient to multiply each equation by the factor $r_{m,n}^{1/4}$.
Let us also denote
\begin{equation}\label{sa.5}\begin{split}
q^+_{m,n}=m^{1/2}r_{m,n}^{1/4}r_{m-1,n}^{1/4},& \qquad
q^-_{m,n}=n^{1/2}r_{m,n}^{1/4}r_{m,n-1}^{1/4};\\
p_{m,n}(\L)&=\z_{m,n}(\L)r^{1/2}_{m,n}.\end{split}
\end{equation}
The equation \eqref{sa.4b} and the similar equation coming from
the condition
\eqref{sa.3} at $x=-1$ yield
\begin{equation}\label{sa.6}\begin{split}
q^+_{m+1,n}C_{m+1,n}^++\frac{2\mu_+p_{m,n}(\L)}{1-e^{-4\z_{m,n}(\L)}}
\left(C_{m,n}^+ -C_{m,n}^-e^{-2\z_{m,n}(\L)}\right)
+q^+_{m,n}C_{m-1,n}^+=0,& \\
q^-_{m,n+1}C_{m,n+1}^-
+\frac{2\mu_-p_{m,n}(\L)}{1-e^{-4\z_{m,n}(\L)}} \left(C_{m,n}^-
-C_{m,n}^+e^{-2\z_{m,n}(\L)}\right) +q^-_{m,n}C_{m,n-1}^-
=0.&\end{split}
\end{equation}
Denote by $\CR=\CR(\L)$ the infinite matrix which corresponds to
this system. In view of \eqref{sa.4a}, we consider $\CR$ as an
operator in the space
\begin{equation*}
 \GG=\ell^2(\N_0^2;\C^2).
\end{equation*}

Removing in \eqref{sa.6} the exponentially small terms, we come to
a simpler system
\begin{equation}\label{sa.7}
q^+_{m+1,n}C_{m+1,n}^+ + 2\mu_+p_{m,n}(\L)C_{m,n}^+
+q^+_{m,n}C_{m-1,n}^+ =0;
\end{equation}
\begin{equation}\label{sa.8}
q^-_{m,n+1}C_{m+1,n}^- + 2\mu_-p_{m,n}(\L)C_{m,n}^-
+q^-_{m,n}C_{m,n-1}^- =0.
\end{equation}
Let $\CR'=\CR'(\L)$ stand for the matrix which corresponds to the
system \eqref{sa.7} -- \eqref{sa.8}, and also for the operator in
$\GG$ generated by this matrix. The operator $\CR'$  decomposes
into an infinite family of simpler operators. First of all, the
equations \eqref{sa.7} (for $C_{m,n}^+$) and \eqref{sa.8} (for
$C_{m,n}^-$) are mutually independent. Further, fix any
$n\in\N_0$. The equations in \eqref{sa.7} which correspond to the
chosen value of $n$ form a linear system in $\ell^2(\N_0)$ with
the Jacobi matrix
\begin{equation*}
\CJ_n^+(\L)=\CJ\left(\{2\mu_+p_{m,n}(\L)\},\{q^+_{m,n}\}\right).
\end{equation*}

In the same way, the equations in \eqref{sa.8}, which correspond
to the chosen value of $m$, form a linear system in $\ell^2(\N_0)$
with the Jacobi matrix
\begin{equation*}
\CJ_m^-(\L)=\CJ\left(\{2\mu_- p_{m,n}(\L)\},\{q^-_{m,n}\}\right).
\end{equation*}
The above reasoning shows that
\begin{equation}\label{sa.11}
\CR'(\L)={\sum_n}^\oplus \CJ_n^+(\L)\oplus{\sum_m}^\oplus \CJ_m^-(\L).
\end{equation}

The original operator $\CR$ can be written as
\begin{equation}\label{sa.11p}
\CR(\L)=\CR'(\L)+\CN(\L)
\end{equation}
where $\CN=\CN(\L)$ is a block-diagonal matrix with $2\times 2$-blocks:
\begin{equation}\label{sa.11n}
 N_{m,n}(\L)=
\frac{2p_{m,n}(\L)e^{-2\z_{m,n}(\L)}}{1-e^{-4\z_{m,n}(\L)}}\begin{pmatrix}
\mu_+e^{-2\z_{m,n}(\L)}&-\mu_+\\ -\mu_-&
e^{-2\z_{m,n}(\L)}&\end{pmatrix}.
\end{equation}
The last two equations elucidate the structure of the matrix
$\CR(\L)$. \vskip0.2cm In the rest of the section we take
$\L=i\tau\in i\R$. We will show that each term on the right-hand
side of \eqref{sa.11} is an invertible operator in $\ell^2$ and
that the norms of $\|\CJ_k^\pm(i\tau)^{-1}\|$ are uniformly
bounded. For this purpose, we note that
\begin{equation*}
(2\mu_+)^{-1}\im\CJ_n^+(i\tau)=(2\mu_-)^{-1}\im\CJ_m^-(i\tau)=
\diag\{p_{m,n}(i\tau)\}.
\end{equation*}
We have $p_{m,n}(i\tau)=\sqrt{r_{m,n}^2-ir_{m,n}\tau}=X+iY$ where
\begin{equation*}
2Y^2=(r_{m,n}^4+r_{m,n}^2\tau^2)^{1/2}-r_{m,n}^2=\frac{r_{m,n}\tau^2}
{(r_{m,n}^2+\tau^2)^{1/2}+r_{m,n}}\ge 2c^2|\tau|.
\end{equation*}
The last inequality, with some constant $c>0$, is valid for
$|\tau|\ge\tau_0$ and for any $m,n\ge 0$; we have taken into
account that $r_{m,n}\ge r_{0,0}=(\nu_+^2+\nu_-^2)/2$.

By the well known estimate for the operators with sign-defined
imaginary part, see e.g. Theorem IV.4.1 in \cite{GK}, this implies
that
\begin{equation*}
\|\CJ_k^\pm(i\tau)^{-1}\|\le\left(c\sqrt{|\tau|}\right)^{-1},
\qquad\forall k\in\N_0,
\end{equation*}
and therefore
\begin{equation}\label{sa.13}
\|\CR'(\tau)^{-1}\|=\sup\|\CJ_k^\pm(\tau)^{-1}\|\le
\left(c\sqrt{|\tau|}\right)^{-1}.
\end{equation}

\vskip0.2cm

The norms of the blocks $N_{m,n}(i\tau)$ in \eqref{sa.11n} are
controlled by $|p_{m,n}(i\tau)|e^{-2\g_{m,n}(i\tau)}$ and hence,
are bounded uniformly in $m,n\in\N_0$. Therefore,
\[\|\CR(\pm i\tau)-\CR'(\pm i\tau)\|\le C=C(\tau_0).\]
Choosing $|\tau|$ large enough, we conclude from \eqref{sa.13} that
\begin{equation*}
\|\CR(\pm i\tau)-\CR'(\pm i\tau)\|<\|\CR'(\pm i\tau)^{-1}\|^{-1}.
\end{equation*}
It follows that the operator $\CR(\pm i\tau)$ has bounded inverse,
and, in particular, the system \eqref{sa.6} has only the trivial
solution in $\GG$.

The proof of \thmref{sa.thm} is complete.

\section{Representation of the resolvent $(\BA_\a-\L)^{-1}$}\label{res}
\subsection{Resolvent.}\label{r}
In order to prove \thmref{one.t2}, we need a convenient
representation for both resolvents involved in \eqref{goal}. Here
we do this for the operator $\BA_\a$. We derive an analogue of the
formula (6.6) (the {\it basic formula}) in \cite{NS} or,
equivalently, (6.4) in \cite{S3}. However, there is an important
difference between the techniques we employ here and those in
\cite{NS} and \cite{S3}. The main goal in both cited papers was
the direct study of the a.c. spectrum of the operator
corresponding to the case of one oscillator. To achieve this
objective, the behaviour of the resolvent as the spectral
parameter approaches the real axis was carefully studied. What we
do here is to apply scattering theory, and use the already known
results of \cite{NS} and \cite{S3}. This makes our analysis much
easier. We are able to work with the resolvents for a fixed value
of the parameter $\L$. We exclude $\L$ from the notation, unless
to do so would be confusing.

Let a function $F\in\GH$ have the decomposition
$F\sim\{f_{m,n}\}$. For any $\a\ge0$ let us denote
\begin{equation}\label{au0}
U_{\a}=\bigl(\BA_{\a}-\L\bigr)^{-1}F\sim\{u_{\a;m,n}\},\qquad
V\sim\{v_{m,n}\}=\{u_{\a;m,n}-u_{0;m,n}\}.
\end{equation}
In the notation for $v$ we do not
reflect dependence on the parameter $\a$.

\vskip0.2cm

The operator $(\BA_0-\L)^{-1}$ can be written in an explicit form.
The functions $u_{0;m,n}$ are given by the formula \eqref{res.1},
with $\z=\z_{m,n}(\L)$, see \eqref{sa.4}. It follows that
\begin{equation}\label{res.1amn}
u_{0;m,n}(\pm1)=(2\z_{m,n})^{-1}
\int_\R e^{-\z_{m,n}|t\mp1|}f_{m,n}(t)dt;\qquad [u'_{0;m,n}](\pm1)=0.
\end{equation}

\vskip0.2cm Each function $v_{m,n}(x)$ belongs to the space
$\CF_{\z_{m,n}}$. Hence, the equalities \eqref{1c} are satisfied
for it, again with $\z=\z_{m,n}(\L)$. Using these equalities and
taking into account that $[u'_{\a;m,n}](\pm1)=[v'_{m,n}](\pm1)$,
we find from the matching condition in \eqref{sa.3} at the point
$x=1$ that
\begin{equation}\label{res.3}\begin{split}
-\frac{2\z_{m,n}}{1-e^{-4\z_{m,n}}}&\left(v_{m,n}(1)-e^{-2\z_{m,n}}v_{m,n}(-1)
\right)\\=
\frac{\a_+}{\sqrt2}&\left(\sqrt{m+1}u_{\a;m+1,n}(1)+\sqrt{m}u_{\a;m-1,n}(1)
\right).
\end{split}\end{equation}
As in section \ref{pr}, we let $\mu_\pm=\sqrt2/\a_\pm $. Since
$v_{m,n}=u_{\a;m,n}-u_{0;m,n}$, the equation \eqref{res.3} yields
\begin{equation}\label{res.3q}\begin{split}
\sqrt{m+1}u_{\a;m+1,n}(1)+\frac{2\mu_+
\z_{m,n}}{1-e^{-4\z_{m,n}}}\left(u_{\a;m,n}(1)
-e^{-2\z_{m,n}}u_{\a;m,n}(-1)\right)&\\+\sqrt{m}u_{\a;m-1,n}(1)
=\frac{2\mu_+ \z_{m,n}}{1-e^{-4\z_{m,n}}}\left(u_{0;m,n}(1)
-e^{-2\z_{m,n}}u_{0;m,n}(-1)\right).
\end{split}
\end{equation}

The next step is the same normalization as in section \ref{pr}.
Denote
\begin{equation}\label{xz}
X^\pm_{m,n}=r^{-1/4}_{m,n}u_{0;m,n}(\pm1),\qquad
Z^\pm_{m,n}=r^{-1/4}_{m,n}u_{\a;m,n}(\pm1).
\end{equation}
Each function $v_{m,n}$ can be
represented as in \eqref{sa.4y}, with  $C^\pm_{m,n}=Z^\pm_{m,n}-X^\pm_{m,n}$.
We use a shortened notation for the corresponding elements in $\GG$:
\begin{equation*}
X=\{X^+_{m,n},X^-_{m,n}\}, \
Z=\{Z^+_{m,n},Z^-_{m,n}\},\ C=\{C^+_{m,n},C^-_{m,n}\},\qquad m,n\in\N_0.
\end{equation*}

Multiplying each equation in \eqref{res.3q} by $r^{1/4}_{m,n}$ and
writing out the similar equations coming from the matching
conditions at $x=-1$, we reduce the system to the form
\begin{equation}\label{res.3r}\begin{split}
q_{m+1,n}^+ Z^+_{m+1,n}+\frac{2\mu_+ p_{m,n}}{1-e^{-4\z_{m,n}}}
\left(Z^+_{m,n}-e^{-2\z_{m,n}}Z^-_{m,n}\right)+q_{m,n}^+ Z^+_{m-1,n}&\\
=\frac{2\mu_+ p_{m,n}}{1-e^{-4\z_{m,n}}}
\left(X^+_{m,n}-e^{-2\z_{m,n}}X^-_{m,n}\right);&\\
q_{m,n+1}^- Z^-_{m,n+1}+\frac{2\mu_- p_{m,n}}{1-e^{-4\z_{m,n}}}
\left(Z^-_{m,n}-e^{-2\z_{m,n}}Z^+_{m,n}\right)+q_{m,n}^- Z^-_{m-1,n}&\\
=\frac{2\mu_-
p_{m,n}}{1-e^{-4\z_{m,n}}}\left(X^-_{m,n}-e^{-2\z_{m,n}}X^+_{m,n}
\right).
\end{split}
\end{equation}
This is the non-homogeneous counterpart of the system
\eqref{sa.6}. In order to write it more conveniently, we need more
notation. All the operators introduced below depend on $\L$ and we
always assume that $\L\notin\R$.

\vskip0.2cm Define the operator $\CS=\CS(\L):\GH\to\GG$ by
\begin{equation*}
\CS:F\mapsto\left\{\frac{r_{m,n}^{1/4}}{2}\int_\R
e^{-\z_{m,n}|t-1|}f_{m,n}(t)dt,\ \frac{r_{m,n}^{1/4}}{2}\int_\R
e^{-\z_{m,n}|t+1|}f_{m,n}(t)dt\right\}.
\end{equation*}
According to \eqref{res.1amn} and \eqref{xz}, this can be written as
\begin{equation}\label{cs}
\CS:F\mapsto\{p_{m,n}X^+_{m,n},p_{m,n}X^-_{m,n}\}.
\end{equation}
It follows from the Cauchy-Schwartz inequality and \eqref{2r} that
the operator $\CS$ is bounded.

The diagonal operator
\begin{equation*}
\CP=\CP(\L):\{X^+_{m,n},X^-_{m,n}\}
\mapsto\{p_{m,n}X^+_{m,n},p_{m,n}X^-_{m,n}\}
 \end{equation*}
acts in $\GG$ and is unbounded. Its inverse $\CP^{-1}$ is a
bounded operator.

Further, let $ \CM=\CM(\L)$ be the operator generated by the
block-diagonal matrix, $\CM=\diag\{M_{m,n}\}$, where
\begin{equation}\label{3e}
 M_{m,n}=(1-e^{-4\z_{m,n}})^{-1}\begin{pmatrix}\mu_+&
-\mu_+e^{-2\z_{m,n}}\\ -\mu_- e^{-2\z_{m,n}}&\mu_-\end{pmatrix}.
\end{equation}
Evidently, $\CM$ is bounded in $\GG$.

Finally, we let
\begin{equation}\label{t}
\CT=\CT(\L): \{C_{m,n}^+,C_{m,n}^-\}\mapsto\{r_{m,n}^{1/4}(C_{m,n}^+
\varf_{m,n}^+ +C_{m,n}^-\varf_{m,n}^-)\}.
\end{equation}
This is a bounded operator acting from $\GG$ into $\GH$.

\vskip0.2cm The system \eqref{res.3r} can be written in the
operator form
\begin{equation*}
\CR Z=2\CM\CS F,
\end{equation*}
whence $Z=2\CR^{-1}\CM\CS F$. Here $\CR=\CR(\L)$ is the operator
in $\GG$ which corresponds to the left-hand side of the system
\eqref{res.3r}, or, equivalently, of the homogeneous system
\eqref{sa.6}. We also have $X=\CP^{-1}\CS F$, so that
\begin{equation*}
C=Z-X=\left(2\CR^{-1}\CM-\CP^{-1}\right)\CS F.
\end{equation*}
If $C$ is found from this equation, then evidently $U_\a-U_0=\CT C$.

Now it follows from the construction that
\begin{equation}\label{3f}
(\BA_\a-\L)^{-1}-(\BA_0-\L)^{-1}=\CT\left(2\CR^{-1}\CM
-\CP^{-1}\right)\CS.
\end{equation}
This is the desired representation of the resolvent of the operator $\BA_\a$.

\subsection{On the matrix $\CR(\L)$.}\label{matr}
It was shown in section \ref{pr} that for $\tau$ large enough the
matrix $\CR(i\tau)$ has a bounded inverse. This allowed us to
conclude that the operator $\BA_\a-\L$ has a bounded inverse for
all $\L\notin\R$, and hence $\ker\CR(\L)=\{0\}$ for all such $\L$.
So, the operator $\CR(\L)^{-1}$ is well-defined. However, this
does not imply automatically that this operator is bounded in
$\GG$. We now show that this property is a direct consequence of
the representation \eqref{3f}. Indeed, \eqref{3f} implies that
\[ 2\CR^{-1}\CM-\CP^{-1}=(\CT^*\CT)^{-1}\CT^*
\left((\BA_\a-\L)^{-1}-(\BA_0-\L)^{-1}\right)\CS^*(\CS\CS^*)^{-1}.\]
It is easy to show that for $\L\notin\R_+$ the operators $\CM$ and $\CT^*\CT$ and
$\CS\CS^*$ (acting in $\GG$)
have bounded inverses. This yields the desired result.

\section{Representation of the resolvent  $(\BA^\circ_\a-\L)^{-1}$}\label{reso}
Our aim here is to derive an analogue of the representation
\eqref{3f} for the operator $\BA^\circ_\a$. One possible way to
proceed is to use the decompositions \eqref{one.2}, \eqref{one.3}.
However, we prefer another way, one which parallels our argument
in section \ref{res}. The calculations are easier for
$\BA^\circ_\a$ than for $\BA_\a$.

\vskip0.2cm

For the objects related to the operator $\BA^\circ_\a$ we use the notation
\begin{equation}\label{au00}
U^\circ_{\a}=\bigl(\BA^\circ_{\a}-\L\bigr)^{-1}F\sim\{u^\circ_{\a;m,n}\},\qquad
V^\circ\sim\{v^\circ_{m,n}\}=\{u^\circ_{\a;m,n}-u^\circ_{0;m,n}\}.
\end{equation}
An analogue of \eqref{res.1amn} is given by
\begin{equation}\label{u1a}
u^\circ_{0;m,n}(\pm1)=(2\z_{m,n})^{-1}\int_{\R_\pm}
\left(e^{-\z_{m,n}|t\mp1|}-e^{-\z_{m,n}(1\pm t)}\right)f_{m,n}(t)dt.
\end{equation}

Next, we derive an analogue of \eqref{res.3}. Taking \eqref{1ac}
and $[{u^\circ_{0;m,n}}'](\pm1)=0$ into account, we get from the
matching condition at $x=1$:
\begin{equation*}
-\frac{2\z_{m,n}}{1-e^{-2\z_{m,n}}}v^\circ_{m,n}(1)=\frac{\a_+}{\sqrt
2}
\left(\sqrt{m+1}u^\circ_{\a;m+1,n}(1)+\sqrt{m}u^\circ_{\a;m-1,n}(1)\right).
\end{equation*}
Since $v^\circ_{m,n}=u^\circ_{\a;m,n}-u^\circ_{0;m,n}$, we find, taking, as before,
$\mu_\pm=\sqrt2/\a_\pm$:
\begin{eqnarray}\label{3d}
\sqrt{m+1}u^\circ_{\a;m+1,n}(1)&+&
\frac{2\mu_+\z_{m,n}}{1-e^{-2\z_{m,n}}} u^\circ_{\a;m,n}(1)
+\sqrt{m}u^\circ_{\a;m-1,n}(1) \nonumber \\
&=& \frac{2\mu_+\z_{m,n}}{1-e^{-2\z_{m,n}}} u^\circ_{0;m,n}(1).
\end{eqnarray}
This is much simpler than the system \eqref{res.3q}, which, of
course, merely reflects the special structure of the operator
$\BA^\circ_\a$ as given by \eqref{one.3}. \vskip0.2cm

The normalization, as in section \ref{res}, reduces \eqref{3d} and
the similar equations for $x=-1$ to the form
\begin{equation}\label{3m}\begin{split}
&q^+_{m+1,n}Z^{\circ,+}_{m+1,n}+\frac{2\mu_+p_{m,n}}{1-e^{-2\z_{m,n}}}Z^{\circ,+}_{m,n}+
q^+_{m-1,n}Z^{\circ,+}_{m-1,n}=\frac{2\mu_+p_{m,n}}{1-e^{-2\z_{m,n}}}X^{\circ,+}_{m,n},\\
&q^-_{m,n+1}Z^{\circ,-}_{m,n+1}+\frac{2\mu_-p_{m,n}}{1-e^{-2\z_{m,n}}}Z^{\circ,-}_{m,n}+
q^-_{m,n+1}Z^{\circ,-}_{m,n-1}=\frac{2\mu_-p_{m,n}}{1-e^{-2\z_{m,n}}}X^{\circ,-}_{m,n}.
\end{split}\end{equation}
Here
\[ X^{\circ,\pm}_{m,n}=r_{m,n}^{-1/4}u^\circ_{0;m,n}(\pm1),\qquad
Z^{\circ,\pm}_{m,n}=r_{m,n}^{-1/4}u^\circ_{\a;m,n}(\pm1).\] The
coefficients $q^\pm_{m,n},\ p_{m,n}$ are the same as in
\eqref{res.3r}, being defined in \eqref{sa.5}. By \eqref{u1a}, we
have
\begin{equation*}\begin{split}
2p_{m,n}X^{\circ,+}_{m,n}&=\int_{\R_+}\left(e^{-\z_{m,n}|t-1|}-e^{-\z_{m,n}(t+1)}\right)
f_{m,n}(t)dt,\\
2p_{m,n}X^{\circ,-}_{m,n}&=\int_{\R_-}\left(e^{-\z_{m,n}|t+1|}-e^{\z_{m,n}(t-1)}\right)
f_{m,n}(t)dt.
\end{split}\end{equation*}

Now we define analogues of the operators involved in the equality
\eqref{res.3r}. First of all, $\CR^\circ=\CR^\circ(\L)$ is the
operator in $\GG$, defined by the infinite matrix which
corresponds to the left-hand side of \eqref{3m}. The operator
$\CR^\circ$ can be written in the form similar to \eqref{sa.11p}:
\begin{equation}\label{ro}
\CR^\circ(\L)=\CR'(\L)+\CN^\circ(\L)
\end{equation}
where
\begin{equation}\label{ro0}
\CN^\circ(\L)=\diag \{N^\circ_{m,n}\},\qquad N^\circ_{m,n}=
2\frac{p_{m,n}e^{-2\z_{m,n}}}{1-e^{-2\z_{m,n}}}
\begin{pmatrix}\mu_+&0\\0&\mu_-\end{pmatrix}.
\end{equation}
The self-adjointness of the operator $\BA^\circ_\a$ in $\GH$
implies that the operator $\CR^\circ(\L)$ is invertible for any
$\L\notin\R$.

The operator $\CS^\circ=\CS^\circ(\L):\GH\to\GG$ is a bounded
operator defined by
\begin{equation*}
 \CS^\circ:F\sim\{f_{m,n}\}
\mapsto\{p_{m,n}X^{\circ,+}_{m,n},p_{m,n}X^{\circ,-}_{m,n}\};
\end{equation*}
cf. \eqref{cs}.

The operator $\CM^\circ=\CM^\circ(\L)$ is the bounded operator on
$\GG$ of block-diagonal form
\begin{equation}\label{3ec}
\CM^\circ=\diag\{M_{m,n}\},\qquad M_{m,n}=
(1-e^{-2\z_{m,n}})^{-1}\begin{pmatrix}\mu_+&0\\0&\mu_-\end{pmatrix}.
\end{equation}

Finally, let
\begin{equation*}
\CT^\circ=\CT^\circ(\L):\{C^+_{m,n},C^-_{m,n}\}\mapsto
\left\{r_{m,n}^{1/4}(C^+_{m,n}\varf^{\circ,+}_{m,n}+
C^+_{m,n}\varf^{\circ,-}_{m,n})\right\}
\end{equation*}
This is a bounded operator acting from $\GG$ into $\GH$.
\vskip0.2cm
As in section \ref{res}, we can re-write the system \eqref{3m} as
\begin{equation}\label{3fc}
(\BA^\circ_\a-\L)^{-1}-(\BA^\circ_0-\L)^{-1}=
\CT^\circ\left(2(\CR^\circ)^{-1}\CM^\circ
-\CP^{-1}\right)\CS^\circ.
\end{equation}
\vskip0.2cm Note that for any $\L\neq\overline{\L}$ the operator
$\CR^\circ(\L)^{-1}$ is bounded. The proof is the same as for the
operator $\CR^{-1}(\L)$, see section \ref{matr}.

\section{Proof of \thmref{one.t2}.}\label{proof}
\subsection{The case $\a=0$.}\label{00}
Here we show that
\begin{equation}\label{goal0}
(\BA^\circ_0-\L)^{-3} - (\BA_0-\L)^{-3}\in \GS_1.
\end{equation}
Recall that in the notation of \eqref{au0} and \eqref{au00}
\begin{equation*}
(\BA_0-\L)^{-1}F=\{u_{0,m,n}\},\qquad  (\BA^\circ_0-\L)^{-1}F=\{u^\circ_{0,m,n}\}
\end{equation*}
where the functions $u_{0,m,n}, u^\circ_{0,m,n}$ are given by the equations
\eqref{res.1} and \eqref{res.1b} respectively, with $\z=\z_{m,n}$. So,
both operators are diagonal. Denote by $\Phi_{m,n},\
\Phi^\circ_{m,n}$ their components, and let
\begin{equation}\label{pr0.qq}\begin{split}
\BQ=\diag\{Q_{m,n}\}=&(\BA^\circ_0-\L)^{-1}-(\BA_0-\L)^{-1},\\
Q_{m,n}=&\Phi^\circ_{m,n}-\Phi_{m,n}.
\end{split}\end{equation}
According to
\eqref{res.1f}, each $Q_{m,n}$ is a rank one operator:
\begin{equation}\label{p2}
Q_{m,n}: f_{m,n}(x)\mapsto -(2\z_{m,n})^{-1}g_{m,n}(x)\int_\R
g_{m,n}(t)f_{m,n}(t)dt
\end{equation}
where $g_{m,n}(x)=e^{-\z_{m,n}|x|}$.
It follows from \eqref{p2} that
\[\|Q_{m,n}\|=(2|\z_{m,n}|\re\z_{m,n})^{-1}\le Cr_{m,n}^{-1},\qquad C=C(\L).\]

The norms of $\Phi_{m,n}$ and $\Phi^\circ_{m,n}$ can be easily
estimated (actually, $\|\Phi_{m,n}\|$ can be calculated
explicitly, since this is a convolution operator). By the ``Schur
test'', the norm $\|\BK\|$ of an integral operator in $L^2$ with
the kernel $K(x,t)$ can be estimated as
\begin{equation*}
\|\BK\|^2\le \sup_t\int|K(x,t)|dx\,\sup_x\int|K(x,t)|dt.
\end{equation*}
Applying this to the operators \eqref{res.1} and \eqref{res.1b}, we find that
\begin{equation*}
\|\Phi_{m,n}\|,\ \|\Phi^\circ_{m,n}\|\le Cr_{m,n}^{-1}.
\end{equation*}
Furthermore, the components of the operator $(\BA^\circ_0-\L)^{-3}
- (\BA_0-\L)^{-3}$ are
\begin{equation*}
\Phi_{m,n}^2 Q_{m,n}+\Phi_{m,n}Q_{m,n}\Phi^\circ_{m,n}+
Q_{m,n}(\Phi^\circ_{m,n})^2.
\end{equation*}
The norm of this operator does not exceed $3C^3r_{m,n}^{-3}$ and, since its rank
is no greater than $3$,
its trace class norm does not exceed  $9C^3r_{m,n}^{-3}$.
By \eqref{sa.2}, these numbers
form a convergent double series, and hence, \eqref{goal0} is established.

Note that the exponent $3$ in \eqref{goal0} can not be replaced
by $2$.

\subsection{Difference between the right-hand sides in \eqref{3f}, \eqref{3fc}.}
\label{e} To shorten our notation, let us denote
\begin{equation}\label{end.0}
\BH=\CT\left(2\CR^{-1}\CM
-\CP^{-1}\right)\CS,\qquad
\BH^\circ=\CT^\circ\left((2\CR^\circ)^{-1}\CM^\circ
-\CP^{-1}\right)\CS^\circ.
\end{equation}
Here we show that
\begin{equation}\label{end.1}
\Psi:=\BH^\circ-\BH=\in\GS_1.
\end{equation}
Since all the inverse operators appearing in \eqref{end.0} are bounded,
 we only need to check that
\begin{gather*}
\CT^\circ-\CT,\ \CS^\circ-\CS,\ \CM^\circ-\CM,\
\CR^\circ-\CR\in\GS_1.
\end{gather*}
Here each operator has a block-diagonal structure, with $(2\!\times\!2)$-blocks,
and it is sufficient to estimate the operator norm of each block and to verify
that the corresponding series converge.
\vskip0.2cm
For the operators $\CT^\circ-\CT$ and $\CS^\circ-\CS$ the result
immediately follows from the
definitions of the operators involved and the estimate \eqref{2c}. For the operator
$\CM^\circ-\CM$ the result is evident from the comparison of \eqref{3e} and \eqref{3ec}.
Finally, for $\CR^\circ-\CR$ the result is implied by \eqref{sa.11p} and
\eqref{ro}, if we take into account evident estimates
of the norms of blocks $N_{m,n}$
in \eqref{sa.11n}
and $N^\circ_{m,n}$ in \eqref{ro0}.

\subsection{End of the proof}\label{end}
Unfortunately, the desired inclusion \eqref{goal} is not implied by \eqref{goal0}
and \eqref{end.1} automatically and we need an extra argument in order to finalize the
proof.

Let us denote $\BG=(\BA_0-\L)^{-1}$.
Using also the notation $\BQ$, $\Psi$ as in \eqref{pr0.qq}, \eqref{end.1}, we can
re-write the equalities \eqref{3f} and \eqref{3fc}  as
\begin{equation*}
(\BA_\a-\L)^{-1}=\BG+\BH,\qquad
(\BA^\circ_\a-\L)^{-1}=\BG+\BH+\BQ+\Psi.
\end{equation*}
We already know that $\Psi\in\GS_1$ and
$(\BG+\BQ)^3-\BG^3\in\GS_1$. Therefore, the following equality is
satisfied modulo a trace class correction:
\begin{gather*} (\BA^\circ_\a-\L)^{-1}-(\BA_\a-\L)^{-1}=
(\BG+\BQ+\BH+\Psi)^3-(\BG+\BH)^3\\= (\BG+\BQ+\BH)^3-(\BG+\BH)^3\
(\mod\GS_1)\\
=\left((\BG+\BQ)^2-\BG^2\right)\BH+(\BG+\BQ)\BH(\BG+\BQ)-\BG\BH\BG\\
+\BH\left((\BG+\BQ)^2-\BG^2\right)+\BQ\BH^2+\BH\BQ\BH+\BH^2\BQ\
(\mod\GS_1).\end{gather*}

Removing the parentheses, we come to the sum where each term
involves one of the products $\BQ\BH$, $\BH\BQ$, $\BQ\BG\BH$,
$\BH\BG\BQ$. Taking into account the structure of the operator
$\BH$, we see that it is sufficient for us to prove that the
operators
\begin{equation}\label{fin}
 \BQ\CT,\ \CS\BQ,\ \BQ(\BA_0-\L)^{-1}\CT ,\
 \CS(\BA_0-\L)^{-1}\BQ
\end{equation}
are trace class. All these operators have block-diagonal form,
with the blocks given by explicit formulae implied by the
corresponding definitions. For instance, according to \eqref{t}
the operator $\BQ\CT$ transforms the number sequence
$\{C_{m,n}^+,C_{m,n}^-\}$ into the sequence of functions
$\{w_{m,n}\}$ where
\begin{equation*}
w_{m,n}(x)=-\frac{r_{m,n}^{1/4}}{2\z_{m,n}}g_{m,n}(x)\int_\R
\left(C_{m,n}^+ \varf_{m,n}^+(t)+,C_{m,n}^- \varf_{m,n}^-(t)\right)g_{m,n}(t)dt.
\end{equation*}
An elementary calculation shows that the integral here is of order
$O(e^{-\g_{m,n}})$. This happens because the function
$g_{m,n}(t)$, see \eqref{p2}, is concentrated around the point
$t=0$, while $\varf_{m,n}^\pm(t)$ is concentrated around $t=\pm1$,
and all the three functions decay exponentially when $t$ moves
away from the corresponding center. Clearly, this estimate implies
that $\BQ\CT\in\GS_1$. The proofs for the other operators in
\eqref{fin} are similar.

The proof of is complete.

\section{Acknowledgments}
The work on the paper started in the Summer of 2004 when one of
the authors (M.S.) was a guest of the School of Mathematics,
Cardiff University. M.S. takes this opportunity to express his
gratitude to the University for its hospitality and to the EPSRC
for financial support under grant GR/T01556.

\end{document}